\documentclass[12pt]{amsart}
\usepackage[margin=1.2in]{geometry}
\usepackage[alphabetic,nobysame]{amsrefs}
\usepackage{amssymb,mathrsfs,amsmath,amscd,multirow,amsthm}
\usepackage{color}
\usepackage[all]{xy}
\usepackage{bigdelim}
\usepackage[inline,shortlabels]{enumitem}%
  \setlist{topsep=3pt, itemsep=0pt,partopsep=0pt}
\usepackage[colorlinks]{hyperref}

\newtheorem{thm}{Theorem}[section]
\newtheorem{cor}[thm]{Corollary}
\newtheorem{prop}[thm]{Proposition}
\newtheorem{ques}{Question}

\theoremstyle{definition}
\newtheorem{defn}[thm]{Definition}
\newtheorem{examp}[thm]{Example}

\theoremstyle{remark}
\newtheorem{rmk}[thm]{Remark}

\newtheorem*{rmk*}{Remark}

\theoremstyle{plain}

\date{\today}
\begin{document}
\title[Varieties with nef anticanonical divisor]{
The structure of varieties with nef anticanonical divisor
in positive characteristic
}

\author{Jiarong Fan}
\email{
fjr24@mail.ustc.edu.cn}
\address{School of Mathematical Science, University of Science and Technology of China, Hefei 230026, P.R.China.}

\author{Lei Zhang}
\email{zhlei18@ustc.edu.cn}
\address{School of Mathematical Science, University of Science and Technology of China, Hefei 230026, P.R.China.}
\begin{abstract}
Let $X$ be a regular projective variety with nef anticanonical divisor over an algebraically closed
field of characteristic $p>0$, and let $f:X\to A$ be a surjective morphism to an abelian variety. We prove that
there exists an isogeny $A' \to A$ such that, the base change $f': X\times_AA' \to A'$ admits a Zariski locally trivial structure.
If moreover $-K_X\equiv D$ for some semiample
$\mathbb Q$-divisor $D$, then $-K_X\sim_{\mathbb Q}D$,
and there exists an isogeny $A' \to A$ such that, the base change $ X\times_AA' \cong X_0\times A'$ where $X_0$ is the fiber over the identity $0\in A'$.
\end{abstract}

\maketitle

\section{Introduction}

Throughout this paper, we work over an algebraically closed field $k$
of characteristic $p>0$.

The geometry of projective varieties with
nef anti-canonical divisors is closely related to
the structure of their Albanese morphisms.
Over $\mathbb{C}$, a smooth projective variety $X$ with  $-K_X$  nef has a beautiful structure:
\begin{itemize}
\item[(i)] the Albanese morphism $a_X:X \to A$ is locally trivial; that is, for any $z\in A$ there exists an analytic neighborhood $U$ of $z$ such that $X\times_AU \cong F\times U$,
where $F$ is a fiber of $a_X$ (conjectured in \cite[Section~4]{DPS96} in the K\"ahler case and proved in the projective case \cite[Theorem~1.2]{Cao19}); 
\item[(ii)] if, moreover, $K_X\equiv 0$, then $K_X\sim_{\mathbb{Q}} 0$,
and there exists an isogeny  $\pi:A' \to A$ such that $X\times_A A' \cong F\times A'$ \cite[Theorem~8.2 and Proposition~8.3]{Kaw85}.
\end{itemize}

In positive characteristic, a typical pathology is that for a fibration $f:X\to Y$ (a projective morphism of normal varieties such that $f_*\mathcal{O}_X = \mathcal{O}_Y$), even if $X$ is smooth, the geometric generic fiber $X_{\bar{\eta}}$ can be singular or even non-reduced. It seems not reasonable to expect analogous results hold in positive characteristic. To formulate a reasonable expectation, we first recall the related results in positive characteristic. Let $X$ be a smooth projective variety with $-K_X$ nef. Then the following results have been established.
\begin{itemize}
\item
The Albanese morphism $a_X$ is surjective.
Let $X \xrightarrow{g} Y \xrightarrow{h} A$ be its Stein factorization. Then $h$ is a finite purely inseparable morphism. 
If either $g$ or $h$ is separable, then $h$ is an isomorphism; see \cite[Theorem~1.3]{EP23}.
\item 
If the geometric generic fiber $X_{\bar{\eta}}$ of $a_X$ is strongly $F$-regular, then all the closed fibers of $a_X$  are isomorphic to each other; if moreover $-K_X$ is numerically equivalent to a semiample $\mathbb Q$-divisor, then there exists an isogeny $A' \to A$ such that
$X\times_A A' \cong X_0\times A'$, where $X_0$ is the fiber of $a_X$ over the identity element; see \cite[Theorems~5.7 and~5.8]{Eji23}. Remark that if $X$ is globally $F$-split, then the geometric generic fiber of of $a_X$ is automatically strongly $F$-regular
by \cite[Theorems~1.2 and~1.3]{Eji19} and \cite[Theorem~6.3]{PZ26}, and $X$ admits a refined decomposition structure; see \cite[Theorem~1.13]{PZ26}, \cite[Theorem~1.5]{Eji23} and \cite[Theorem~1.7]{Zha26}.
\item
If $\dim A=1$ and $-K_X$ is numerically equivalent to a semiample $\mathbb Q$-divisor, then
there exists an isogeny  $\pi:A'\to A$ such that 
\mbox{$X\times_A A' \cong F\times A'$}
for a closed fiber $F$ of $a_X$;  see \cite[Theorem~6.1]{Eji23}. In case $\dim X =2$ and $K_X \equiv 0$, this result should trace back to Bombieri and Mumford \cite{BM77, BM76}.

\item If $a_X:X\to A$ has relative dimension one,
then it is a fibration \cite[Theorem~1.5]{CWZ26b}.
If, in addition,
$K_X\equiv 0$, then
there exists an isogeny  \mbox{$\pi:A'\to A$}
such that $(X\times_A A')^{\nu}_{\mathrm{red}} \cong C\times A'$, where $C$ is an elliptic curve or $C=\mathbb{P}^1$ \cite[Theorems 5.2--5.4]{CWZ26a}.
If instead $K_X\not\equiv 0$, then there exists an isogeny 
$\pi:A'\to A$ such that $(X\times_A A')^{\nu}_{\mathrm{red}}$ is a $\mathbb P^1$-bundle over $A'$ \cite[Theorems~1.1 and 1.2]{GLZ25}. These varieties often can be realized by the quotient of certain foliations. 
\end{itemize}

According to the above progress, we propose the following two questions. We expect both of them have positive answers and unify all the above results (Example \ref{ex:unify} is a typical example).

\begin{ques}\label{ques:1}
Let $X$ be a regular projective variety over $k$ with $K_X\equiv 0$ or, more generally, $-K_X\equiv D$ for some  semiample $\mathbb Q$-divisor $D$. Is $-K_X$ semiample?
\end{ques}

\begin{ques}\label{ques:2}
Let $X$ be a regular projective variety over $k$ with $-K_X$  nef. Does the Albanese morphism $a_X:X\to A$  admit a certain splitting structure?
\end{ques}

In this paper, we answer both questions
affirmatively in full generality.

\begin{thm}\label{thm:abundance}
Let $X$ be a normal projective variety over $k$.
Let $\Delta$ be an effective $\mathbb Q$-Weil divisor on $X$ such that $(X,\Delta)$ is
strongly $F$-regular.
Suppose that  
\begin{itemize}
    \item $K_X+\Delta$ has Cartier index not divisible by $p$, and
    \item $-(K_X+\Delta)\equiv D$ for a semiample $\mathbb Q$-Cartier $\mathbb Q$-divisor $D$.
\end{itemize}

Then $-(K_X+\Delta)\sim_{\mathbb Q}D$.
Moreover, for every surjective morphism $f:X\to A$ 
to an abelian variety,
    there exists an isogeny $\pi:A'\to A$ such that
    $$ (X,\Delta)\times_A A'
 \cong (X_0,\Delta|_{X_0})\times_k A'
    $$
 as $A'$-schemes.
 Here, $X_0$ is the fiber of $f$ over the origin $0\in A$, and
 $\Delta|_{X_0}$ is regarded as an effective $\mathbb Q$-almost Cartier divisor.
\end{thm}

Under the weaker assumption
that $-(K_X+\Delta)$ is nef,
we prove that
every such morphism $f:X\to A$ admits a similar splitting structure locally.

\begin{thm}\label{thm:split-av}
Let $X$ be a normal projective variety over $k$, and let $\Delta$ be an
effective $\mathbb Q$-Weil divisor such that $(X,\Delta)$ is strongly
$F$-regular. Suppose that $-m(K_X+\Delta)$ is a nef Cartier divisor for
some integer $m>0$ not divisible by $p$. Let $f\colon X\to A$ be a
surjective morphism to an abelian variety, 
and denote by $X_0$ the fiber of $f$ over the origin $0\in A$.

Then  there exist an isogeny $\pi:A'\to A$
and a Zariski open cover $A'=\bigcup U_i$
with $U_i$-isomorphisms
    $$ (X,\Delta)\times_A U_i
 \cong (X_0,\Delta|_{X_0})\times_k U_i.
    $$
\end{thm}

We also  show that,
in the setting of 
 Theorem \ref{thm:split-av},
 the Albanese morphism of $X$ is a fibration.

\begin{cor}\label{cor:albanese-fibration}
Let $X$ be a normal projective variety over $k$, and let $\Delta$ be an
effective $\mathbb Q$-Weil divisor such that $(X,\Delta)$ is strongly
$F$-regular. Suppose that $-m(K_X+\Delta)$ is a nef Cartier divisor for
some integer $m>0$ not divisible by $p$. 

Then the Albanese morphism $f:X\to A$
is a fibration.
\end{cor}

\textbf{Strategy of the proof}.
We briefly explain the proofs of Theorems~\ref{thm:abundance} and~\ref{thm:split-av}.
By Propositions \ref{prop:drepo-to-drestr} and \ref{prop:torsion-and-splitting} (see also \cite[Theorem~1.5]{Zha26}), we have been done once we can find an $f$-ample Cartier divisor $B$ on $X$
such that 
$f_*\mathcal O_X(nB)$ is numerically flat for every $n>0$. This has been proved by Ejiri \cite{Eji23} under the additional assumption that 
the geometric generic fiber of $f$, together with the induced boundary, is strongly $F$-regular. \textbf{In fact this is enough to treat the situation without the $F$-regularity condition on the fiber.} The precise result we shall use is 
\begin{thm}[{\cite[Theorem~4.2]{Eji23}}]\label{thm:Eji}
Let $(X,\Delta)$ be a strongly
$F$-regular projective pair over
 $k$. Let $f\colon X\to A$ be a
surjective morphism to an abelian variety.
Suppose that
\begin{itemize}
    \item[$(*)$] $-m(K_X+\Delta)$ is a nef Cartier divisor for
some integer $m>0$ not divisible by $p$, and
\item[$(**)$] $(X_{\bar{\eta}},\Delta|_{X_{\bar{\eta}}})$
is strongly $F$-regular,
where $X_{\bar{\eta}}$ is the geometric generic fiber of $f$.
\end{itemize}
Let $B$ be an $f$-ample Cartier divisor on $X$.
Then there exists an integer $l_0>0$ such that, if setting $r=\mathrm{rank} f_*\mathcal O_X(l_0B)$, then for every $n>0$,
$$f_*\mathcal O_X(nl_0B)\langle-\frac{n}{r}\det f_*\mathcal O_X(l_0B)\rangle$$ is a numerically flat $\mathbb{Q}$-twisted bundle.  In particular, for the $f$-ample divisor $C = rl_0B- f^*\det f_*\mathcal O_X(l_0B)$ and integer any $n>0$, $f_*\mathcal O_X(nC)$ is  a numerically flat
vector bundle on $A$. 
\end{thm}
We shall apply Theorem \ref{thm:Eji} to 
the projection $p_2: \widetilde{X}= X\times \widehat A\to \widehat A$ and the polarization $\widetilde{B}=p_1^*B \otimes (f\times\mathrm{id}_{\widehat A})^* \mathcal{P}_A$, where $B$ is an ample line bundle on $X$, $\widehat A =  \mathrm{Pic}^0(A)$ and $\mathcal{P}_A$ is the Poincar\'{e} line bundle on $A\times \widehat A$. Theorem \ref{thm:Eji} tells  that there exists a $\mathbb{Q}$-divisor $D$ on $\widehat{A}$ such that for sufficiently divisible $n$, $p_{2*}(\widetilde{B}^{\otimes n})\langle-nD\rangle$ is numerically flat. From this, by the proof of \cite[Theorem~5.8]{Eji23} (or Proposition \ref{prop:torsion-and-splitting}), we conclude that for any $\alpha \in \widehat{A}$ we have a  polarized isomorphism $h_{\alpha}: (X,B) \cong (X, B\otimes f^*\mathcal{P}_{\alpha})$. This has been quite close to our final conclusion. For example, when $K_X \equiv 0$, this roughly gives an action of $\widehat{A}$ on $X$ (up to an isogeny), and one can then use the treatment of \cite[Section 3]{BM77}.

\begin{examp}\label{ex:unify}
To  illustrate Theorem~\ref{thm:abundance},
let's  consider the following example, which is a special case of \cite[Example~14.1]{PZ26}.
Let $E$ be a supersingular elliptic curve over $k$, and let
$ A=E\times E$.
Let $\mathrm{Fr}_A^{(1)}\colon A\to A^{(1)}$ be the relative
Frobenius.
Then
$\ker(\mathrm{Fr}_A^{(1)})\cong\alpha_p^{2}$, where
$\alpha_p=\operatorname{Spec}(k[\epsilon]/(\epsilon^p))$ is 
the subgroup of $p$-nilpotents in $\mathbb G_a$.
 Every line
$\ell\in
    \mathbb P\bigl(\operatorname{Lie}(\alpha_p^{2})\bigr)$
determines a subgroup scheme
$G_\ell$ of $\alpha_p^{2}$, and together these subgroup schemes
form a tautological family $\mathcal G\subseteq A\times\mathbb P^1$.
Let
\[
    \pi\colon A\times\mathbb P^1
    \longrightarrow
    X:=(A\times\mathbb P^1)/\mathcal G
\]
be the quotient.
Then $X$ is  smooth and 
is an abelian scheme over $\mathbb P^1$,
with structure morphism
 $g\colon X\to\mathbb P^1$
 induced by the projection
$p_2\colon A\times \mathbb P^1\to \mathbb P^1$.
The composition
$\mathrm{Fr}_A^{(1)}\circ p_1\colon A\times \mathbb P^1\to A^{(1)}$ induces 
the Albanese morphism  $f\colon X \to A^{(1)}$. 
Moreover, 
$$\omega_X\cong g^*\mathcal O_{\mathbb P^1}(p-3).
$$
Thus $-K_X$ is semiample for $p=2,3$.
As predicted by Theorem~\ref{thm:abundance},
the Albanese morphism $f\colon X\to A^{(1)}$ splits into a
product after base change along an isogeny.

We now construct the splitting
explicitly.
Since translation on the
first factor of $A\times\mathbb P^1$ commutes with the action of $\mathcal G$, it induces an action of $A$ on $X$ and thereby gives
an $A$-isomorphism
\[
X_0\times A \overset{\sim}{\longrightarrow}X\times_{A^{(1)}}A,
\quad (x,a)\longmapsto (a\cdot x,a),
\]
where $X_0$ is the  fiber of $f$ over the origin
$0\in A^{(1)}$.
Diagrammatically, we have
$$\centerline{\xymatrix@C=1.8cm{ 
&  X\cong (A\times \mathbb P^1)/\mathcal G\ar[d]^-{f}
& X \times_{A^{(1)}}A\cong X_0 \times A\ar[l]\ar[d]
\\
&  A^{(1)} & A\ar[l]_-{\mathrm{Fr}_A^{(1)}}
}}
$$

For $p=3$,
 \cite[Proposition~14.3]{PZ26} states that
\textbf{if
$Z\to X$ is a finite étale cover  and
$Y\to Z$ is a torsor under a finite infinitesimal $k$-group scheme with $Y$ integral,
then the Albanese morphism of $Y$ does not split}.
This
shows that even when $X$ is smooth,
the fiber  $X_0$ need not be reduced.
We can describe $X_0$ 
explicitly from the viewpoint of 
foliations as in
\cite[Section~6.2]{CWZ26a}.
Indeed, since $\mathcal G$ has height one,
$X$ coincides with
the quotient of $A\times \mathbb P^1$ by 
the rank-one
foliation $\mathcal F$ induced by the Lie algebra of $\mathcal G$.
Let $\alpha$ and $\beta$ be two invariant vector
fields on $A$
such that $\alpha(0)$ and $\beta(0)$
form a basis of
$\mathrm{Lie}(\alpha_p^2)$.
Let $[u:v]$ be the corresponding
 homogeneous coordinates 
on $\mathbb P^1\cong \mathbb{P}(\mathrm{Lie}(\alpha_p^2))$.
By definition, $\mathrm{Lie}(G_{[u:v]})$ is spanned by
 $u\alpha(0)+v\beta(0)$. 
Hence  $\mathcal F$  is 
generated by $\alpha+(v/u)\beta$ on $\{ u\ne 0\}$ and by $(u/v)\alpha+\beta$ on $\{ v\ne 0\}$.
Choose linear coordinates $x,y$ on
\mbox{$\ker (\mathrm{Fr}_A^{(1)})\cong \alpha_p^2$} whose differentials at
the origin form the basis dual to $\alpha(0),\beta(0)$.
Write
 $t=v/u$ and $s=u/v$ for the
coordinates on the two standard affine charts of $\mathbb P^1$.
Set $R=k[t,x,y]/(x^p,y^p)$. Then over $U=\{u\ne0\}$, we have
$$X_0\times_{\mathbb P^1} U\cong
\mathrm{Spec}\bigl(
\mathrm{Ann}_R (\mathcal F|_{\ker(\mathrm{Fr}_A^{(1)})\times U})
\bigr)
\cong\mathrm{Spec}\bigl(k[t,z_U]/(z_U^p)\bigr),
\quad z_U=y-tx.
$$
Similarly,
over $V=\{v\ne0\}$, we have
$$X_0\times_{\mathbb P^1} V\cong
\mathrm{Spec}\bigl(k[s,z_V]/(z_V^p)\bigr),
\quad z_V=sy-x.
$$

\end{examp}

\subsection{Further remarks and questions.}
Varieties over imperfect fields are also crucial even in the classification of varieties over algebraically closed fields. For instance, to study the Iitaka fibration $f: X \to Y$ of a good minimal model $X$ over an algebraically closed field, the generic fiber $X_{\eta}$ of $f$ is such a variety over $k(\eta)$ with $K_{X_{\eta}} \equiv 0$ and the geometric generic fiber $X_{\bar{\eta}}$ being possibly quite singular. 

\begin{ques} Let $X$ be a regular projective variety  over an imperfect field $K=H^0(\mathcal{O}_X)$, with $-K_X$ being nef.

(1) If $-K_X$ is numerically semiample, is it semiample?

(2) If $K$ is separably closed, which implies that the Albanese variety has a rational point, is the  Albanese morphism $X\to A$ locally trivial in the sense of 
Theorem \ref{thm:split-av}?
\end{ques}

\begin{rmk}
The semi-ampleness of a divisor and $F$-regularity are preserved under separable field extensions. For these questions, we can pass to the separable closure of $K$. 
However, we cannot apply Ejiri's Theorem \ref{thm:Eji}  directly. In our situation, tracing back to the proof of \cite[Theorem~4.2]{Eji23}, what is  essentially needed is that the fiber product 
$$X^{(r)} = X\times_K X \times_K \cdots \times_K X$$
is $F$-regular. But this does not hold if $X$ is not geometrically $F$-regular over $K$.
\end{rmk}

\medskip
{\small \noindent\textit{\textbf{AI statement.}}
We learned from  GPT-6-Astra the key idea: to consider the projection
$X\times \widehat A\to \widehat A$. The model originally suggested studying the map $\tilde{f}: X \times A \to A,\,(x,a) \mapsto f(x)-a$. 
We have reformulated the statements and reorganized the proof with assistance from AI tools (GPT and Kimi). }
\medskip

{\small \noindent\textit{\textbf{Acknowledgments.}}
This research is partially supported by CAS Project for Young Scientists in Basic Research (No. YSBR-032) and NSFC (No. 12122116 and No. 12471495). The authors would like to thank Yang Zhang for his useful comment and  reminding us the stimulating \cite[Example~14.1]{PZ26}.
}

\section{Notation and conventions}
\begin{itemize}
\item 
We work over an algebraically closed field $k$ with characteristic $p>0$.
    \item
    By a $k$-scheme, we mean a separated scheme of finite type over $k$.
    A variety is an
    integral $k$-scheme.
 A projective morphism $f : X \to S$ between normal  varieties is called a fibration if
    $f_*\mathcal{O}_X = \mathcal{O}_S$.

\item
For a  coherent sheaf $\mathcal E$ of rank $r$
on a normal variety $X$, we write $\det\mathcal E=(\bigwedge^r\mathcal E)^{\vee\vee}$,
where $(-)^{\vee\vee}$ denotes the double dual.
An open subset $U$ of a $k$-scheme $Y$ is called big if
$\operatorname{codim}_Y(Y\setminus U)\geq 2$.

    \item
    We use $\equiv$, $\sim$ and $\sim_{\mathbb Q}$ for numerical, linear, and
    $\mathbb Q$-linear equivalence, respectively.

    \item
    For a variety $X$, we denote by
    $F_X^e\colon X\to X$ and
    $\mathrm{Fr}_X^{(e)}\colon X\to X^{(e)}:=X\times_{k,F_k^e}k$
    the $e$-fold absolute Frobenius morphism of $X$ and the $e$-fold
    relative Frobenius morphism of $X$ over $k$, respectively.

    \item Let $X$ be a normal  variety over $k$,
and let $\Delta $ be an effective $\mathbb Q$-Weil divisor
on $X$.
We say that $(X, \Delta)$ is a strongly F-regular pair if for every effective Weil divisor $D$ on $X$, there exists an $e\in \mathbb{Z}_{>0}$ such that the natural homomorphism
\[
\mathcal{O}_X \to F_{X*}^e \mathcal{O}_X \hookrightarrow F_{X*}^e\mathcal{O}_X(\lceil (p^e-1) \Delta  \rceil+D)
\] splits locally as an $\mathcal{O}_X$-module homomorphism (\cite[Definition 2.1]{HW02}).
For $e\in \mathbb{Z}_{>0}$, 
applying $\mathcal{H}om(-,\mathcal{O}_X)$ to the natural homomorphism
$$\mathcal{O}_X \longrightarrow F^e_{X*}\mathcal{O}_X 
\hookrightarrow 
F^e_{X*}\mathcal{O}_X (\lfloor (p^e-1)\Delta\rfloor)
$$
gives, by Grothendieck duality, the trace map 
\[
{\mathrm Tr}^{e}_{X,\Delta}: F^e_{X*}\mathcal O_X(
\lceil (1-p^e)(K_X+\Delta)\rceil)\to \mathcal{O}_X.
\] 
For a Cartier divisor $B$ on $X$, we denote by 
${\mathrm Tr}^{e}_{X,\Delta}(B)$ the homomorphism obtained by
tensoring ${\mathrm Tr}^{e}_{X,\Delta}$ with $\mathcal O_X(B)$.

\item
Let $X$ be an equi-dimensional $k$-scheme satisfying
$S_2$ and $G_1$.
Here, $S_2$ (resp.\ $G_1$) is Serre's second condition
(resp.\ the condition that it is Gorenstein in codimension $1$).
Let $\mathcal K$ be the sheaf of total quotient rings on $X$.
An AC divisor (or almost Cartier divisor) on $X$
is a reflexive coherent subsheaf of $\mathcal K$ which is
invertible on a big open subset of $X$
(\cite{Har94}*{p.~301}, \cite{MS12}*{Definition~2.1}).

Let $D$ be an AC divisor on $X$,
and let $\mathcal O_X(D)$ denote the coherent sheaf defining $D$.
We say that $D$ is effective if
$\mathcal O_X\subseteq\mathcal O_X(D)$.
The set $\mathrm{WSh}(X)$ of AC divisors on $X$ forms
naturally an additive group
(\cite{Har94}*{Corollary~2.6}).
A $\mathbb Q$-AC divisor is an element of
$\mathrm{WSh}(X)\otimes_{\mathbb Z}\mathbb Q$.
A $\mathbb Q$-AC divisor is called effective if it can be written
as $D\otimes \lambda$, where $D$ is an effective AC divisor
and $\lambda\in\mathbb Q_{\geq 0}$.

\item
Let $S$ be a $k$-scheme, and 
let $f_i:X_i\to S$ $(i=1,2)$ be projective morphisms,
where $X_i$ are equi-dimensional $k$-schemes satisfying $S_2$ and $G_1$.
Let $B_i$ be an $f_i$-ample Cartier divisor and
$\Delta_i$ an effective $\mathbb Q$-AC divisor on $X_i$.
We write
\[
(X_1,B_1,\Delta_1)\cong_S(X_2,B_2,\Delta_2)
\]
for  polarized isomorphism over $S$ (\cite[Definition~7.2]{PZ26});
that is, there exist an $S$-isomorphism $\varphi:X_1\to X_2$
and a line bundle $\mathcal M$ on $S$ such that
\[
\varphi^*\mathcal O_{X_2}(B_2)
\cong \mathcal O_{X_1}(B_1)\otimes f_1^*\mathcal M,
\qquad
\varphi^*\Delta_2=\Delta_1.
\]
    \item
    For an abelian variety $A$, let $\widehat A=\operatorname{Pic}^0(A)$,
    let $t_a$ denote translation by $a$, and let $[n]_A$ denote
   the multiplication-by-$n$ morphism on $A$.
    We denote by $\mathcal P_A$ the normalized Poincar\'e line bundle 
    on 
    $A\times \widehat A$
   and by $\mathcal P_{\alpha}$ the line bundle on $A$
   corresponding to 
    $\alpha \in\widehat A$.
   
\end{itemize}

\section{Numerical flatness and isotriviality}\label{sec:split}

\subsection{Numerically flat vector bundles on abelian
varieties}

\begin{defn}
Let $\mathcal E$ be a vector bundle on an abelian variety $A$
over $k$.
\begin{enumerate}
    \item We say that $\mathcal E$ is numerically flat if both $\mathcal E$
and $\mathcal E^\vee$ are nef.
\item We say that $\mathcal E$ is homogeneous if 
$t_a^*\mathcal E\cong\mathcal E$
for every $a\in A$.
\end{enumerate}
\end{defn}

\begin{thm}\label{thm:flat-vb-av}
Let $\mathcal E$ be a vector bundle on an abelian variety $A$ over  $k$. Then the following are equivalent:
\begin{enumerate}
    \item $\mathcal E$ is numerically flat;
    \item $\mathcal E$ is  an iterated extension of algebraically trivial line bundles;
    \item $\mathcal E$ is homogeneous;
    \item $\mathcal E \simeq \bigoplus_i \mathcal U_i \otimes \mathcal L_i$, where
    each  $\mathcal L_i\in \mathrm{Pic}^0(A)$
    and each $\mathcal U_i$ is an iterated extension of $\mathcal O_A$;
   \item  there exists an isogeny of abelian varieties $\pi:A'\to A$
   such that $\pi^*\mathcal E\cong\bigoplus \mathcal P_i$,
   where each $\mathcal P_i\in \mathrm{Pic}^0(A')$.
\end{enumerate}
\end{thm}
\begin{proof}
$(1)\Rightarrow(2)$ is established in
\cite[\S6]{Lan12}.
 $(2)\Rightarrow(4)$
 can be proved by induction on the rank.
 $(4)\Rightarrow(1)$ is trivial.
 $(3)\Rightarrow(4)$ 
 has been proved in \cite[Theorem~2.3]{Miy73}.
 $(4)\Rightarrow(3)$ has been proved in
\cite[Theorem~4.17]{Muk78}.
 $(4)\Rightarrow(5)$ follows from the fact that the map $[p]_A^*:H^1(A,\mathcal O_A)\to H^1(A,\mathcal O_A)$
 induced by $[p]_A:A\to A$ is zero.
 $(5)\Rightarrow(1)$ follows from the fact that
  numerical flatness descends under isogenies.
\end{proof}

\subsection{Auxiliary Positivity} The following positivity follows easily from a combination of  Fujita vanishing \cite[Theorem 1.5]{Kee03} and Castelnuovo-Mumford regularity \cite[Lemma 2.5.1]{FGA05}, which is a standard strategy.

\begin{prop}\label{prop:nef-direct-images}
Let $(X,\Delta)$ be a strongly
$F$-regular projective pair over
 $k$.
Suppose that $-m(K_X+\Delta)$ is a nef Cartier divisor for some integer $m>0$ not divisible by $p$.
Let $f:X\to A$ be a flat surjective morphism to an abelian variety,
and let $B$ be a nef and $f$-ample Cartier divisor on $X$.

Then there exists an integer $n_0>0$ such that
for every $n\ge n_0$ and every nef Cartier divisor $N$ on $X$,
$$
f_*\mathcal O_X(nB+N)
$$
is a nef vector bundle.
\end{prop}

\begin{proof}
We can choose  an integer $n_0>0$ such that the following properties hold.
\begin{enumerate}
  
\item For every nef Cartier divisor $N$ on $X$,
every  $e>0$ with $(p^e-1)(K_X+\Delta)$ Cartier,
and every $n\ge n_0$,
the homomorphism
$$f_*\mathrm{Tr}^{e}_{X,\Delta}(nB+N):
F^{e}_{A*}f_*\mathcal O_X
\bigl((1-p^{e})(K_X+\Delta)+p^{e}(nB+N)\bigr)
\longrightarrow f_*\mathcal O_X(nB+N)
$$
is surjective.
 This
 follows from
\cite[Proposition~3.1]{Eji23}.
  \item For every  nef Cartier divisor $N$ on $X$
  and every $n\ge n_0$,
  $f_*\mathcal O_X(nB+N)$
  is locally free.
This follows from Keeler's relative Fujita vanishing \cite[Theorem 1.5]{Kee03}
and cohomology and base change.
\end{enumerate}

Fix an integer $n\ge n_0$ and a nef Cartier divisor $N$.
Let $\ell $ be a positive integer not divisible by $p$,
let
$$f_{\ell}:X_{\ell}:=X\times_{A,[\ell]}A\to A
$$
 be the base change of $f$
by $[\ell]:A\to A$,
and let $\rho_{\ell}:X_{\ell}\to X$ be the projection.
Set $B_{\ell}:=\rho^*_{\ell}B$,
$N_{\ell}:=\rho^*_{\ell}N$,
  $\Delta_{\ell}:=\rho^*_{\ell}\Delta$,
  and
$$
\mathcal G_{e,\ell}:=
F^e_{A*}f_{\ell*}\mathcal O_{X_\ell}
\bigl(
 (1-p^e)(K_{X_\ell}+\Delta_\ell)+p^e(nB_\ell+N_\ell)\bigr).
$$
Since $[\ell]$ is étale, 
for each
$e>0$ with $(p^e-1)(K_X+\Delta)$ Cartier,
the homomorphism 
$$[\ell]^*f_*\mathrm{Tr}^{e}_{X,\Delta}(nB+N):\mathcal G_{e,\ell}\longrightarrow [\ell]^*f_*\mathcal O_X(nB+N)$$
is also surjective.
Let $H$ be a  very ample divisor  on $A$ such that
$B+f^*H$ is ample. Then $B_{\ell}+f_{\ell}^*H$ is also ample.
By Keeler's relative Fujita vanishing \cite[Theorem 1.5]{Kee03}
and
Castelnuovo-Mumford regularity \cite[Lemma 2.5.1]{FGA05}, 
for $e\gg0$,
$$\mathcal G_{e,\ell}\otimes\mathcal O_A(n(\dim A+1)H)$$
is globally generated.
The surjectivity of $[\ell]^*f_*\mathrm{Tr}^{e}_{X,\Delta}(nB+N)$
implies that
$$[\ell]^*f_*\mathcal O_X(nB+N)\otimes\mathcal O_A(n(\dim A+1)H)$$
is also globally generated, and therefore nef.
Letting $\ell$ tend to infinity shows that $f_*\mathcal O_X(nB+N)$ is nef.
\end{proof}

\subsection{Isotriviality}

We show how to deduce isotriviality from the numerical flatness of direct images.

\begin{prop}\label{prop:drepo-to-drestr}
Let $(X,\Delta)$ be a strongly
$F$-regular projective pair over
 $k$.
Suppose that $-m(K_X+\Delta)$ is a nef Cartier divisor for some integer $m>0$ not divisible by $p$.
Let $f: X\to A$ be a surjective morphism to an abelian variety,
and denote by $X_0$ the fiber of $f$ over the origin $0$.

Let
$B$ be an $f$-ample Cartier divisor on $X$
such that for every $n>0$, $f_*\mathcal O_X(nB)$
is a numerically flat vector bundle.
Then there exist
    an integer $d>0$, 
    an isogeny $A'\to A$
and an open cover $A'=\bigcup U_i$
with polarized $U_i$-isomorphisms
$$(X,dB,\Delta)\times_A U_i\cong_{U_i} (X_0,dB|_{X_0},\Delta|_{X_0})\times_k U_i.
$$
In particular, $f$ is flat.
\end{prop}

\begin{proof}
 After replacing $B$ by a sufficiently large multiple, we may assume that
for every $n>0$, the natural homomorphism
$\mathrm{Sym}^n (f_*\mathcal{O}_X(B)) \to
    f_*\mathcal{O}_X(nB)$ is surjective.
Set
$$\mathcal{K}_n:=\ker(\mathrm{Sym}^n (f_*\mathcal{O}_X(B)) \to
    f_*\mathcal{O}_X(nB)).$$
There exists an integer $n_0>0$ such that $\bigoplus_{n>0} \mathcal{K}_n$ is generated in degrees at most $n_0$.

By Theorem \ref{thm:flat-vb-av},
there exists an isogeny $\pi:A'\to A$ such that for every $0<n\le n_0$,
$\pi^*f_*\mathcal{O}_X(nB)\cong \bigoplus_j \mathcal  L_{nj}$,
where each $\mathcal L_{nj}\in \mathrm{Pic}^0(A')$.
We may moreover assume that the subgroup 
of $\mathrm{Pic}^0(A')$
generated by the
$\mathcal L_{1j}$ is torsion-free.
Note that for any two line bundles $\mathcal{L},\mathcal{L}'\in \mathrm{Pic}^0(A')$,
$$
    \mathrm{Hom}(\mathcal{L}, \mathcal{L}') = 
    \begin{cases}
        k, & \mathcal{L} \cong \mathcal{L}', \\
        0, & \mathcal{L} \not\cong \mathcal{L}'.
    \end{cases}
$$
For each $0<n \le n_0  $, since both 
$\mathrm{Sym}^n (\pi^*f_*\mathcal{O}_X(B))$ and $\pi^*f_*\mathcal{O}_X(nB)$
are direct sums of line bundles in $\mathrm{Pic}^0(A')$,
the natural homomorphism 
$$
    \mathrm{Sym}^n (\pi^*f_*\mathcal{O}_X(B)) \rightarrow \pi^*f_*\mathcal{O}_X(nB)
$$
 is represented by a matrix with entries in $k$.
Hence  $\pi^*\mathcal K_n$
is a direct sum of certain summands of 
$ \mathrm{Sym}^n (\pi^*f_*\mathcal{O}_X(B))\cong \bigoplus (\mathcal L_{1j_1}\otimes\cdots \otimes \mathcal L_{1j_n})$.
Consequently,
for every $n>0$,
$\pi^*\mathcal{K}_n$ and $\pi^*f_*\mathcal{O}_X(nB)$
are  direct sums of 
line bundles of the form 
$\mathcal L_{1j_1}\otimes\cdots \otimes \mathcal L_{1j_n}$.

Set $\mathcal R:=\mathcal O_{A'}\oplus
\bigoplus_{n>0}\pi^*f_*\mathcal O_X(nB)$.
Take an open cover $A'=\bigcup_i U_i$ trivializing all
$\mathcal L_{1j}$.
Then 
$$\mathcal R|_{U_i}
\cong
(\mathcal R\otimes k(0))\otimes_k\mathcal O_{U_i}$$
as graded $\mathcal O_{U_i}$-algebras.
Hence
$$(X,B)\times_A U_i\cong_{U_i}(X_0,B|_{X_0})\times_k U_i.
$$

Finally, we show that these polarized isomorphisms can be chosen to preserve
the boundary.
Fix an integer $m>0$ such that
$-m(K_X+\Delta)$
is  Cartier.
 Since $t_a^*\mathcal R\cong \mathcal R$ as 
graded $\mathcal O_{A'}$-algebras
for every $a\in A'$, $\pi^*f_*\mathcal O_X(nB+mK_X)$ is a homogeneous vector bundle on $A'$ for every $n>0$.
Hence $f_*\mathcal O_X(nB+mK_X)$ is numerically flat for every $n>0$.
By relative Serre vanishing and Proposition \ref{prop:nef-direct-images}, 
for $n\gg0$,
the natural homomorphism
$$f_*\mathcal O_X(nB-m(K_X+\Delta))\otimes
 f_*\mathcal O_X(nB+mK_X)\longrightarrow
 f_*\mathcal O_X(2nB-m\Delta)
$$
is surjective with nef source.
Hence, the target
$$f_*\mathcal O_X(2nB-m\Delta)$$
is  nef,
and therefore numerically flat, since it is
a subsheaf of the numerically flat vector bundle $f_*\mathcal O_X(2nB)$.
Consequently, 
after replacing $\pi$ by a further isogeny,
we may assume that, for a fixed $n\gg0$, 
$\pi^*f_*\mathcal O_X(2nB-m\Delta)$ is a direct sum
of line bundles in $\mathrm{Pic}^0(A')$.
The inclusion
\[
\pi^*f_*\mathcal O_X(2nB-m\Delta)
\hookrightarrow \pi^*f_*\mathcal O_X(2nB)
\]
is then represented by constant matrices with respect to the
decompositions into algebraically trivial line bundles.
Thus the 
polarized isomorphisms above identify the pullback of
$\mathcal O_X(-m\Delta)$ with the pullback of 
$\mathcal O_{X_0}(-m\Delta|_{X_0})$ as ideal sheaves.
Hence they preserve the boundary $\Delta$.
\end{proof}

When $-(K_X+\Delta)$ is semiample, we can obtain a global
splitting structure, which refines \cite[Theorem~5.8]{Eji23} by using almost the same
 argument.

\begin{prop}\label{prop:torsion-and-splitting}
Let $X$ be a normal projective variety over $k$, and let
$\Delta$ be an effective $\mathbb Q$-Weil divisor on $X$ such that
$(X,\Delta)$ is strongly $F$-regular.
Let $f:X\to A$ be a surjective morphism to an abelian variety,
and denote by $X_0$ the fiber over the origin $0\in A$.
Let $B$ be an $f$-ample Cartier divisor on $X$.
Suppose that
\begin{enumerate}
\item $K_X+\Delta$ has Cartier index not divisible by $p$;
\item $-m(K_X+\Delta)\sim D+f^*P$ for some integer $m>0$,
a basepoint-free Cartier divisor $D$ on $X$, and a Cartier divisor
$P$ on $A$ with $\mathcal O_A(P)\in\mathrm{Pic}^0(A)$;
\item $f_*\mathcal O_X(nB)$ is numerically flat for every $n>0$.
\end{enumerate}

Then $P\sim_{\mathbb Q}0$.
Moreover,
 there exist an integer $d>0$
 and an
 isogeny
$\pi:A'\to A$ with a polarized $A'$-isomorphism
\[
(X,dB,\Delta)\times_A A'\cong_{A'} (X_0,dB|_{X_0},\Delta|_{X_0})\times_k A'.
\]
\end{prop}

\begin{proof}
Since $f_*\mathcal O_X(B)$ is an iterated extension of line bundles
in $\mathrm{Pic}^0(A)$, after replacing
$B$ by $B-f^*Q$ for some divisor $Q$ on $A$ with $\mathcal O_A(Q)\in\mathrm{Pic}^0(A)$, we may  assume that $B$ is effective.
Let $a$ be the Cartier index of $K_X+\Delta$,
and set
\[
I:=\left\{i\in\mathbb Z  \mid 0\le i<m,\ a\mid i\right\}.
\]
Choose an integer $d>1$ not divisible by $p$ such that
$(X,\Delta+\frac{1}{d}B)$ is strongly $F$-regular.
Replacing $B$ and $d$ by $nB$ and $nd$ for a sufficiently large
integer $n$ not divisible by $p$, we can assume that the following
properties hold.
\begin{enumerate}
\item[(P1)] For every $i\in I$ and every $e>0$ with $ad\mid(p^e-1)$,
the homomorphism
\begin{align*}
&f_*\mathrm{Tr}_{X,\Delta+\frac{1}{d}B}^e(B-i(K_X+\Delta)):\\
&F^e_{A*}f_*\mathcal O_X((1-p^e)(K_X+\Delta+\frac{1}{d}B)+p^e(B-i(K_X+\Delta)))
\longrightarrow f_*\mathcal O_X(B-i(K_X+\Delta))
\end{align*}
 is surjective.
 This
 follows from
\cite[Proposition~3.1]{Eji23}.
\item[(P2)] For every $i\in I$ and every $v\ge1$,
the natural homomorphism
\[
f_*\mathcal O_X(B-i(K_X+\Delta))\otimes
(f_*\mathcal O_X(B))^{\otimes(v-1)}
\longrightarrow f_*\mathcal O_X(vB-i(K_X+\Delta))
\]
is surjective.
This follows from
 Keeler's relative Fujita vanishing \cite[Theorem~1.5]{Kee03} and the relative Castelnuovo-Mumford regularity \cite[Example~1.8.24]{Laz04}.
\item[(P3)] For every $i\in I$, $v\ge1$ and $q\ge0$,
there is a surjective homomorphism
\[
\bigoplus f_*\mathcal O_X(vB-i(K_X+\Delta))
\longrightarrow f_*\mathcal O_X(qD+vB-i(K_X+\Delta)).
\]
This follows
from \cite[Lemma~5.9]{Eji23}.
\item[(P4)] For every $i\in I$, 
$$\mathcal G_i:=f_*\mathcal O_X(B-i(K_X+\Delta))$$
 is a nef vector bundle. This follows from Proposition \ref{prop:nef-direct-images}.
\end{enumerate}
Combining (P2) and (P3),
we obtain the following:
\begin{enumerate}
    \item[(P5)]  For every $i\in I$, $v\ge1$ and $q\ge0$,
there is a surjective homomorphism
\[
\Phi_{q,v,i}:\bigoplus (\mathcal G_i\otimes \mathcal G_0^{\otimes (v-1)})
\longrightarrow f_*\mathcal O_X(qD+vB-i(K_X+\Delta)).
\] 
\end{enumerate}
Fix $e\ge m+1$ such that $p^e-1=ads$ for some $s\in \mathbb Z_{>0}$.
For each $i\in I$,
write $(i+1)p^e-1=mq_i+r_i$, with $q_i\in\mathbb Z_{>0}$ and $r_i\in [0,m-1]$.
Then $r_i\in I$, and 
\begin{align*}
&(1-p^e)(K_X+\Delta+\frac{1}{d}B)+p^e(B-i(K_X+\Delta))\\
=&(p^e-as)B-((i+1)p^e-1)(K_X+\Delta)\\
\sim & (p^e-as)B-r_i(K_X+\Delta)+q_i(D+f^*P).
\end{align*}
Applying (P5) with $q=q_i$, $v=p^e-as$ and $i=r_i$,
we obtain a surjective homomorphism
\begin{align*}
\Psi_i:&\bigoplus \bigl(\mathcal G_{r_i}\otimes \mathcal G_0^{\otimes (p^e-as-1)}\otimes \mathcal O_A(q_iP)\bigl)\\
&\longrightarrow f_*\mathcal O_X((p^e-as)B-r_i(K_X+\Delta)+q_i(D+f^*P))\\
& ~\cong 
f_*\mathcal O_X((1-p^e)(K_X+\Delta+\frac{1}{d}B)+p^e(B-i(K_X+\Delta))).
\end{align*}
Combining $F^e_{A*}\Psi_i$ with the surjection in (P1),
we obtain a surjection
\[
\Psi_i':\bigoplus F^e_{A*}\bigl(\mathcal G_{r_i}\otimes \mathcal G_0^{\otimes (p^e-as-1)}\otimes \mathcal O_A(q_iP)\bigl)
\longrightarrow\mathcal G_i.
\]

For $i\in I$, set
\[
 E(i):=\{\mathcal L\in\mathrm{Pic}^0(A)\mid
 \mathrm{Hom}(\mathcal G_i,\mathcal L)\ne0\}.
\]
Since $\mathcal G_i$ is nef, $E(i)$ is finite by \cite[Lemma 4.5]{Eji22}.
Set
\[
M(1):=\bigcup_{i\in I}E(i)=\{\mathcal Q_1,\dots,\mathcal Q_\lambda\},
\]
and, for every $n\ge1$, set
\[
M(n):=\left\{
\bigotimes_{j=1}^{\lambda}\mathcal Q_j^{\otimes n_j}
\ \middle|\ n_j\in\mathbb Z_{\ge0},\
\sum_{j=1}^{\lambda}n_j=n
\right\}.
\]
Since $B$ is effective,
$\mathcal O_A\in E(0)$.
Hence $M(n)\subseteq M(n+1)$.

For $\mathcal L\in E(i)$, composing $\Psi_i'$ with a nonzero
homomorphism $\mathcal G_i\to\mathcal L$
and using Grothendieck duality,
we obtain a nonzero homomorphism
\[
\Psi_i'':\mathcal G_{r_i}\otimes\mathcal G_0^{\otimes(p^e-as-1)}
\otimes\mathcal O_A(q_iP)
\longrightarrow\mathcal L^{\otimes p^e}.
\]
Since $\mathcal G_0^{\otimes(p^e-as-1)}$ admits a filtration
$$0=\mathcal E_0\subset \mathcal E_1\subset \cdots \subset \mathcal E_g=
\mathcal G_0^{\otimes(p^e-as-1)}
$$
with $\mathcal E_j/\mathcal E_{j-1}\in M(p^e-as-1)$,
$\Psi''_i$ induces a nonzero homomorphism
$$\mathcal G_{r_i}\otimes(\mathcal E_j/\mathcal E_{j-1})
\otimes\mathcal O_A(q_iP)
\longrightarrow\mathcal L^{\otimes p^e}
$$
for some $j$.
It follows that
\[
\mathcal L^{\otimes p^e}\otimes\mathcal O_A(-q_iP)\in M(p^e-as).
\]

By the pigeonhole principle, it follows from the above relation that
for every  line bundle $\mathcal M\in M(n)$
with  $n>p^e\cdot|M(1)|$,
there exists an $i\in I$ such that
\begin{equation}\tag{$*$}
\mathcal M\otimes\mathcal O_A(-q_iP)
\in M(n-as)\subseteq M(n).
\end{equation}
Since $M(n)$ is finite, iterating $(*)$ gives
 $P\sim_{\mathbb Q}0$.
Moreover, $(*)$ implies that
 $\bigcup_{n\ge1}M(n)$ is finite.
 In particular, each line bundle in $E(0)$ is torsion.
Hence the isogeny $\pi:A'\to A$ in the proof of Proposition \ref{prop:drepo-to-drestr}
can be chosen so that for every $n>0$,
$\pi^*f_*\mathcal O_X(nB)$ is a trivial vector bundle,
and therefore
$$(X,B,\Delta)\times_A A'\cong (X_0,B|_{X_0},\Delta|_{X_0})\times_k A'.$$
\end{proof}

\section{Proofs of the main results}

In view of Propositions~\ref{prop:drepo-to-drestr} and~\ref{prop:torsion-and-splitting}, to prove
Theorems~\ref{thm:abundance} and~\ref{thm:split-av}, it suffices to establish the following theorem, which has been proved by Ejiri  under the condition that the geometric generic fiber is strongly $F$-regular; see Theorem~\ref{thm:Eji}. It is surprising that Theorem~\ref{thm:Eji} implies that an analogous result holds without the $F$-regularity condition on the fiber.

\begin{thm}\label{thm:numerically-flat-direct-images}
Let $(X,\Delta)$ be a strongly $F$-regular projective pair
over $k$.
Suppose that $-m(K_X+\Delta)$ is a nef Cartier divisor
for some integer $m>0$ not divisible by $p$.
Let $f:X\to A$ be a surjective morphism to an abelian variety.

Then there exists an $f$-ample  Cartier divisor $B$
on $X$ such that, for every $n\ge0$,
\[
f_*\mathcal O_X(nB)
\]
is a numerically flat vector bundle on $A$.
\end{thm}

\begin{proof}
Let $p_X:X\times \widehat A\to X$ and 
$p_{\widehat A}:X\times \widehat A\to \widehat A$ be the projections.
Choose a very ample divisor $B_0$ on $X$.
Let
$\mathcal P_A$ be the Poincar\'e bundle on
$A\times\widehat A$,
and set
$$
\mathcal O_{X\times \widehat A}(\widetilde B)
:=p_X^*\mathcal O_X(B_0)\otimes
(f\times\operatorname{id}_{\widehat A})^*\mathcal P_A.
$$
Then $\widetilde B$ is $p_{\widehat A}$-ample.
Since 
$(X,\Delta)$ is strongly $F$-regular and
$k$ is algebraically closed, the geometric generic
fiber of
$(X\times\widehat A,p_X^*\Delta)\to\widehat A$,
with the induced boundary,
is strongly $F$-regular.
Applying Theorem~\ref{thm:Eji}  to the morphism $p_{\widehat A}:\widetilde X:=X\times \widehat A\to \widehat A$ and the polarization $\widetilde B$,
we obtain that there exist
an integer $d>0$ and a divisor $D$ on $\widehat{A}$ such that for every $n>0$ the sheaf 
$$(p_{\widehat A})_*\mathcal{O}_{\widetilde X}(n(d\widetilde B -p_{\widehat A}^*D))$$
is a numerically flat vector bundle. Then applying  Proposition~\ref{prop:drepo-to-drestr}, there exist
an isogeny
$\sigma:A'\to\widehat A$, and an open neighborhood
$U'\subseteq A'$ of the origin
with a polarized $U'$-isomorphism
$$(X\times \widehat A,d\widetilde B)\times_{\widehat A} U'\cong_{U'} (X,dB_0)\times_k U'.
$$
In  other words,
we obtain a family of automorphisms
\[
h:U'\longrightarrow\operatorname{Aut}(X),
\quad u\mapsto h_u
\]
such that, for every $u\in U'(k)$,
\begin{align*}\tag{$*$}
    h_u^*\mathcal O_X(dB_0)
\cong
\mathcal O_X(dB_0)\otimes f^*\mathcal P_{d\cdot\sigma(u)}.
\end{align*}
Replacing $h_u$ by $h_0^{-1}\circ h_u$, we may assume
that $h_0=\operatorname{id}_X$.

Set $r:=\mathrm{rank}(f_*\mathcal O_X(dB_0))$ and 
$$\mathcal O_X(B):=\mathcal O_X(rdB_0)\otimes f^*\bigl(\det f_*\mathcal O_X(dB_0)\bigr)^{-1}.$$
We shall show that for every $n\ge0$, $f_*\mathcal O_X(nB)$
is a homogeneous vector bundle, and hence it is numerically flat by Theorem \ref{thm:flat-vb-av}.
 
Since the morphism 
$$\tilde{f}: X\times U'\to A,(x,u)\mapsto f(h_u(x))-f(x)$$
sends $X\times \{0_{A'}\}$ to $0_A$,
by the rigidity lemma \cite[p.~43]{Mumford-AV}, after shrinking $U'$ if necessary,
we obtain a factorization $\tilde{f} = \pi\circ p_{U'}: X\times U'\xrightarrow{p_{U'}} U' \xrightarrow{\pi} A$, in other words
for each closed point $u\in U'$ we have 
$f(h_u(x))-f(x) = \pi(u)$, which gives
the following commutative diagram 
$$\centerline{\xymatrix@C=2cm{ 
&X\ar[d]_{f} &  X\ar[l]_-{h_u}\ar[d]^{f}\\
& A & A\ar[l]_{t_{\pi(u)}}.
}}
$$
Then, for every $n\ge0$ and every $u\in U'(k)$, we have
\[
\begin{aligned}
t_{\pi(u)}^*f_*\mathcal O_X(nB)
&\cong
t_{\pi(u)}^*f_*\mathcal O_X(nrdB_0)
\otimes t_{\pi(u)}^*\bigl(\det f_*\mathcal O_X(dB_0)\bigr)^{-n}\\
&\cong
f_*h_u^*\mathcal O_X(nrdB_0)
\otimes
\bigl(\det f_*h_u^*\mathcal O_X(dB_0)\bigr)^{-n}\\
&\cong
\bigl(f_*\mathcal O_X(nrdB_0)
\otimes\mathcal P_{nrd\cdot\sigma(u)}\bigr)
\otimes
\bigl(\det f_*\mathcal O_X(dB_0)
\otimes\mathcal P_{rd\cdot\sigma(u)}\bigr)^{-n}\\
&\cong f_*\mathcal O_X(nB).
\end{aligned}
\]
Here the first isomorphism follows from the projection formula,
the second follows from flat base change, and the third  follows from $(*)$.

Thus, the desired homogeneity will follow once we show that $\pi(U')$ generates $A$.
Note that 
 $\pi$ is quasi-finite.
 Indeed,
 for any closed points $u,v\in U'$ with $\pi(u)=\pi(v)$, $(*)$ gives
\begin{align*}
    \mathcal{P}_{rd\cdot\sigma(u)}\cong &
\bigl(\det  f_*h_u^*\mathcal O_X(dB_0)\bigr)\otimes
\bigl(\det f_*\mathcal O_X(dB_0)  \bigr)^{-1}\\
\cong &
\bigl(\det t_{\pi(u)}^* f_*\mathcal O_X(dB_0)\bigr)\otimes
\bigl(\det f_*\mathcal O_X(dB_0)  \bigr)^{-1} \\
\cong &
\bigl(\det t_{\pi(v)}^* f_*\mathcal O_X(dB_0)\bigr)\otimes
\bigl(\det f_*\mathcal O_X(dB_0)  \bigr)^{-1}\\
\cong &
\bigl(\det  f_*h_v^*\mathcal O_X(dB_0)\bigr)\otimes
\bigl(\det f_*\mathcal O_X(dB_0)  \bigr)^{-1}\\
\cong &
\mathcal P_{rd\cdot \sigma(v)}.
\end{align*}
Hence $u-v$ belongs to  $\ker([rd]_{\widehat A}\circ\sigma)$, which is finite.
This shows that  $\pi$ is quasi-finite,
and therefore dominant.
Then by Chevalley's theorem,
$\pi(U')$
contains a nonempty open subset $U\subset A$.
For every $a\in A(k)$,
since
$U\cap (a+U)\ne \varnothing$,
we can write $a=u_1-u_2$ with $u_1,u_2\in U(k)$.
Hence 
$$t_a^*f_*\mathcal O_X(nB)\cong
(t_{u_1}\circ t_{u_2}^{-1})^*f_*\mathcal O_X(nB)\cong
f_*\mathcal O_X(nB).
$$

\end{proof}

\begin{rmk}
If $-(K_X+\Delta)$ is semiample,
using Proposition~\ref{prop:torsion-and-splitting} in place of Proposition~\ref{prop:drepo-to-drestr} in the above proof,
we can take $U'$ to be the whole $A'$.
Thus we obtain an  action of $A'$  on $X$, and
we can directly construct an $A'$-isomorphism
\[
X_0\times_k A'
\overset{\sim}{\longrightarrow}
X\times_{A,\pi}A',
\quad
(x,a)\longmapsto(h_a(x),a).
\]

However, in Theorem~\ref{thm:abundance}, $-(K_X+\Delta)$ is assumed only
to be numerically equivalent to a semiample $\mathbb Q$-Cartier $\mathbb Q$-divisor;
 its semiampleness is not known a priori.
We therefore first establish Theorem~\ref{thm:numerically-flat-direct-images}, which allows us
to apply Proposition~\ref{prop:torsion-and-splitting} to the Albanese morphism and deduce abundance.
\end{rmk}

\begin{proof}[Proof of Theorem \ref{thm:abundance}]
   Let
    $a_X:X\to \mathrm{Alb}(X)$ be the Albanese morphism of $X$,
 which is surjective by    \cite[Theorem~1.3(a)]{EP23}.
   Set $N:=-(K_X+\Delta)-D\equiv0$.
    By
\cite[Theorem~9.6.3]{FGA05}, there exist an integer $m>0$
and a Cartier divisor $P$ on $\mathrm{Alb}(X)$ with $\mathcal O_{\mathrm{Alb}(X)}(P)\in \mathrm{Pic}^0(\mathrm{Alb}(X))$ such that $mN\sim_{\mathbb Q} a_X^*P$.
Thus,
by Theorem \ref{thm:numerically-flat-direct-images},
we can
apply Proposition \ref{prop:torsion-and-splitting} to $a_X$
to
 obtain $P\sim_{\mathbb Q}0$. Hence $-(K_X+\Delta)\sim_{\mathbb Q}D$.
Moreover, we can therefore apply Proposition~\ref{prop:torsion-and-splitting}
to any surjective morphism $f:X\to A$ to an abelian variety
to obtain the desired splitting structure.

\end{proof}

\begin{proof}[Proof of Theorem \ref{thm:split-av}]
It follows from Theorem \ref{thm:numerically-flat-direct-images}
and Proposition \ref{prop:drepo-to-drestr}.
\end{proof}

\begin{proof}[Proof of Corollary \ref{cor:albanese-fibration}]
By Theorem~\ref{thm:numerically-flat-direct-images},
$f_*\mathcal O_X$ is numerically flat.
 Theorem~\ref{thm:flat-vb-av}  then gives a decomposition
\[
f_*\mathcal O_X\cong
\bigoplus_i \mathcal U_i\otimes \mathcal P_i,
\]
where each $\mathcal U_i$ is an iterated extension of
$\mathcal O_A$ and each $\mathcal P_i\in\operatorname{Pic}^0(A)$.
Each $\mathcal P_i$ is therefore  a subsheaf of $f_*\mathcal O_X$.
The inclusion
$\mathcal P_i\hookrightarrow f_*\mathcal O_X$ induces, by adjunction,
a nonzero homomorphism
$f^*\mathcal P_i\to\mathcal O_X$,
which is an isomorphism since $f^*\mathcal P_i$
is numerically trivial.
The injectivity of
$f^*\colon\operatorname{Pic}^0(A)\to\operatorname{Pic}(X)$
then implies that $\mathcal P_i\cong\mathcal O_A$.
Hence $f_*\mathcal O_X$ is an iterated extension of
$\mathcal O_A$.

In other words, the natural inclusion
$\mathcal O_A\hookrightarrow f_*\mathcal O_X$
extends to a filtration
\[
\mathcal O_A=\mathcal E_1\subset \mathcal E_2\subset\cdots\subset \mathcal E_n
=f_*\mathcal O_X
\]
with $\mathcal E_{j+1}/\mathcal E_j\cong\mathcal O_A$ for every $1\le j<n$.

Suppose  $n>1$. Consider the exact sequence
\[
\tag{*}
0\longrightarrow\mathcal O_A
\longrightarrow \mathcal E_2
\longrightarrow\mathcal O_A
\longrightarrow0.
\]
Pulling back this sequence to $X$, we obtain an exact sequence
\[
0\longrightarrow\mathcal O_X
\longrightarrow f^*\mathcal E_2
\longrightarrow\mathcal O_X
\longrightarrow0,
\]
which splits since the inclusion
$\mathcal O_X\to f^*\mathcal E_2$ admits a retraction given by
the natural composite
\[
f^*\mathcal E_2\longrightarrow f^*f_*\mathcal O_X
\longrightarrow\mathcal O_X.
\]
Then the sequence $(*)$ also splits since
$$f^*\colon H^1(A,\mathcal O_A) \cong \mathrm{Ext}^1_A(\mathcal O_A, \mathcal O_A)\to H^1(X,\mathcal O_X)\cong \mathrm{Ext}^1_X(\mathcal O_X, \mathcal O_X)$$ 
is injective.
Consequently,
$\mathcal O_A^{\oplus2}\subseteq f_*\mathcal O_X$,
contradicting $h^0(A,f_*\mathcal O_X)=1$.
Therefore $f_*\mathcal O_X=\mathcal O_A$.
\end{proof}

\bibliographystyle{plain}

\end{document}